\documentstyle[amssymb,12pt]{amsart}

\newtheorem{cor}{Corollary}
\newtheorem{lem}{Lemma}

\theoremstyle{definition}

\theoremstyle{definition}
\newtheorem{thm}{Theorem}

\theoremstyle{remark}
\newtheorem{rem}{Remark}

\numberwithin{equation}{section}

\begin{document}

\newcommand{\thmref}[1]{Theorem~\ref{#1}}
\newcommand{\secref}[1]{Sect.~\ref{#1}}
\newcommand{\lemref}[1]{Lemma~\ref{#1}}
\newcommand{\propref}[1]{Proposition~\ref{#1}}
\newcommand{\corref}[1]{Corollary~\ref{#1}}
\newcommand{\remref}[1]{Remark~\ref{#1}}
\newcommand{\nc}{\newcommand}
\nc{\on}{\operatorname}
\nc{\ch}{\mbox{ch}}
\nc{\Z}{{\Bbb Z}}
\nc{\C}{{\Bbb C}}
\nc{\pone}{{\Bbb C}{\Bbb P}^1}
\nc{\pa}{\partial}
\nc{\F}{{\cal F}}
\nc{\arr}{\rightarrow}
\nc{\larr}{\longrightarrow}
\nc{\al}{\alpha}
\nc{\ri}{\rangle}
\nc{\lef}{\langle}
\nc{\W}{{\cal W}}
\nc{\la}{\lambda}
\nc{\ep}{\epsilon}
\nc{\su}{\widehat{\goth{sl}}_2}
\nc{\sw}{\goth{sl}}
\nc{\g}{\goth{g}}
\nc{\h}{\goth{h}}
\nc{\n}{\goth{n}}
\nc{\N}{\widehat{\n}}
\nc{\ab}{\goth{a}}
\nc{\G}{\widehat{\g}}
\nc{\De}{\Delta_+}
\nc{\gt}{\widetilde{\g}}
\nc{\Ga}{\Gamma}
\nc{\one}{{\bold 1}}
\nc{\hh}{\widehat{\h}}
\nc{\z}{{\goth Z}}
\nc{\zz}{{\cal Z}}
\nc{\Hh}{{\cal H}}
\nc{\qp}{q^{\frac{k}{2}}}
\nc{\qm}{q^{-\frac{k}{2}}}
\nc{\La}{\Lambda}
\nc{\wt}{\widetilde}
\nc{\qn}{\frac{[m]_q^2}{[2m]_q}}
\nc{\cri}{_{\on{cr}}}
\nc{\k}{h^\vee}
\nc{\sun}{\widehat{\sw}_N}
\nc{\HH}{{\cal H}_q(\sw_N)}
\nc{\ca}{\wt{{\cal A}}_{h,k}(\sw_2)}
\nc{\si}{\sigma}
\nc{\gl}{\widehat{\goth{g}\goth{l}}_2}
\nc{\el}{\ell}
\nc{\s}{t}
\nc{\tN}{\theta_{p^N}}
\nc{\ds}{\displaystyle}
\nc{\Dp}{D_{p^{-1}}}
\nc{\Dq}{D_{q^{-1}}}
\nc{\bi}{\bibitem}

\title{Quantum $\W$--algebras and elliptic algebras}

\author{Boris Feigin}
\address{Landau Institute for Theoretical
Physics, Moscow 117334, Russia and R.I.M.S., Kyoto University, Kyoto 606,
Japan}

\author{Edward Frenkel}\thanks{The research of the second author was
partially supported by NSF grant DMS-9501414}
\address{Department of Mathematics, Harvard University, Cambridge, MA
02138, USA}

\date{}

\maketitle

\section{Introduction.}
\subsection{}
In \cite{FR} N.~Reshetikhin and the second author introduced new Poisson
algebras $\W_q(\g)$, which are $q$--deformations of the classical
$\W$--algebras. The Poisson algebra $\W_q(\g)$ is by definition the center
of the quantized universal enveloping algebra $U_q(\G^L)$ at the
critical level, where $\g^L$ is the Langlands dual Lie algebra to $\g$. It
was shown in \cite{FR} that the Wakimoto realization of
$U_q(\G^L)$ constructed in \cite{AOS} provides a homomorphism
from the center of $U_q(\G^L)$ to a Heisenberg-Poisson algebra
$\Hh_q(\g)$. This homomorphism can be viewed as a free field realization of
$\W_q(\g)$. When $q=1$, it becomes the well-known Miura transformation
\cite{DS}. In \cite{FR} explicit formulas for this free field realization
were given. The structure of these formulas is the same as that of the
formulas for the spectra of transfer-matrices in integrable quantum spin
chains obtained by the Bethe ansatz method \cite{KS}. This is not
surprising given that these spectra can actually be computed using the
center of $U_q(\G^L)$ at the critical level and the Wakimoto
realization. For the Gaudin models, which correspond to the $q=1$ case,
this was explained in detail in \cite{FFR}.

\subsection{}
The Poisson algebra $\W_q(\sw_2)$ is a $q$--deformation of the classical
Virasoro algebra. It has generators $t_n, n\in\Z$, and relations \cite{FR}
\begin{equation}    \label{virrel}
\{ t_n,t_m \} = - h \sum_{l\in\Z} \frac{1-q^l}{1+q^l} t_{n-l} t_{m+l} - h
(q^n-q^{-n}) \delta_{n,-m},
\end{equation}
where $h=\log q$. The Heisenberg-Poisson algebra $\Hh_q(\sw_2)$
has generators $\la_n, n \in \Z$, and relations
\begin{equation}    \label{hrel}
\{ \la_n,\la_m \} = - h \frac{1-q^n}{1+q^n} \delta_{n,-m}.
\end{equation}
Let us form the generating series $$\La(z) = q^{-1/2} \exp \left(
- \sum_{n\in\Z} \la_n z^{-n} \right), \quad \quad t(z) = \sum_{n\in\Z} t_n
z^{-n}.$$ The free field realization homomorphism is given by
\begin{equation}    \label{ff}
t(z) \longrightarrow \Lambda(zq^{1/2}) + \Lambda(zq^{-1/2})^{-1}.
\end{equation}
As shown in \cite{FR}, in the limit $q \arr 1$, the algebra $\W_q(\sw_2)$
becomes isomorphic to the classical Virasoro algebra while formula
\eqref{ff} becomes the Miura transformation.

In \cite{SKAO} J.~Shiraishi, H.~Kubo, H.~Awata, and S.~Odake quantized
formulas \eqref{virrel}, \eqref{hrel}, and \eqref{ff}. They constructed a
non-commutative algebra depending on two parameters $q$ and $p$, such that
when $q=p$ it becomes commutative, and is isomorphic to the Poisson algebra
$\W_q(\sw_2)$. Let us denote this algebra by $\W_{q,p}(\sw_2)$. Shiraishi,
e.a., constructed a free field realization of $\W_{q,p}(\sw_2)$, i.e. a
homomorphism into a Heisenberg algebra $\Hh_{q,p}(\sw_2)$, which has the
same form as \eqref{ff}. They also constructed the screening currents,
i.e. operators acting on the Fock representations of $\Hh_{q,p}(\sw_2)$,
which commute with the action of $\W_{q,p}(\sw_2)$ up to a total
difference. They showed that if one fixes $\beta\in\C$ and sets
$p=q^{1-\beta}$, then in the limit $q \arr 1$ the algebra $\W_{q,p}(\sw_2)$
becomes isomorphic to the Virasoro algebra with central charge
$1-6(1-\beta)^2/\beta$ \cite{SKAO}.

The work of Shiraishi, e.a. \cite{SKAO} was motivated by their
bosonization formula for the Macdonald symmetric functions
\cite{AOS:bos}. The paper \cite{SKAO} reveals a remarkable connection
between the algebra $\W_{q,p}(\sw_2)$ and Macdonald's functions
corresponding to rectangular Young diagrams: those turn out to coincide
with singular vectors of $\W_{q,p}(\sw_2)$ in its bosonic Fock
representations.

\subsection{}
The goal of the present work is to generalize the results of Shiraishi,
e.a., to the case of the $q$--deformed $\W$--algebras, and to point out
some intriguing elliptic structures arising in these algebras. Namely, we
construct an algebra $\W_{q,p}(\sw_N)$ depending on $q$ and $p$, such that
when $q=p$ it becomes isomorphic to the $q$--deformed classical
$\W$--algebra $\W_q(\sw_N)$ from \cite{FR}. We also construct, along the
lines of \cite{FR} and \cite{SKAO}, a free field realization of
$\W_{q,p}(\sw_N)$, which is a deformation of the free field realization
from \cite{FR}, and the screening currents. One can observe many
similarities between the algebra $\W_{q,p}(\sw_N)$ and the ordinary
$\W$--algebra of $\sw_N$ constructed by V.~Fateev and S.~Lukyanov \cite{FL}
(see also \cite{BG}), which can be recovered from $\W_{q,p}(\sw_N)$ in the
limit $q \arr 1$.

The algebra $\W_{q,p}(\sw_N)$ is topologically generated by Fourier
coefficients of currents $T_1(z),\ldots,T_{N-1}(z)$. The free field
realization of $\W_{q,p}(\sw_N)$ is defined by the formula
\begin{multline}    \label{freefield}
\Dp^N - T_1(z) \Dp^{N-1} + T_2(z) \Dp^{N-2} - \ldots + (-1)^{N-1}
T_{N-1}(z) \Dp + (-1)^N = \\ :(\Dp - \La_1(z))(\Dp - \La_2(zp)) \ldots (\Dp
- \La_N(zp^{N-1})):,
\end{multline}
where $\La_i(z), i=1,\ldots,N$, are generating series of a Heisenberg
algebra, and $[\Dp \cdot f](x) = f(xp^{-1})$. In the limit $q \arr 1$ this
formula becomes the normally ordered Miura transformation from \cite{FL}.

The screening currents $S^\pm_i(z)$ are solutions of the difference
equations:
\begin{align*}
D_q S_i^+(z) & = p^{-1} :\La_{i+1}(zp^{i/2}) \La_i(zp^{i/2})^{-1}
S^+_i(z):,\\ D_{p/q} S_i^-(z) & = p^{-1} :\La_{i+1}(zp^{i/2})
\La_i(zp^{i/2})^{-1} S^-_i(z):.
\end{align*}
Using formula \eqref{freefield} one can check that they commute with the
currents $T_i(z)$ up to a total difference. This implies that their
residues acting between bosonic Fock representations commute with the
action of $\W_{q,p}(\sw_N)$.

Using these operators one can construct singular vectors in the Fock
representations of $\W_{q,p}(\sw_2)$. These singular vectors should give
the Macdonald symmetric functions corresponding to Young diagrams with
$N-1$ rectangles as was pointed out in \cite{SKAO}.

\subsection{}
An interesting aspect of the algebras $\W_{q,p}(\sw_N)$ is the appearance
of elliptic functions in their definition and free field realization.

In particular, we show that the series $\La_i(z), i=1,\ldots,N$, satisfy,
in the analytic continuation sense, the following relations:
\begin{equation}
\La_i(z) \La_j(w) = \varphi_N\left( \frac{w}{z} \right) \La_j(w) \La_i(z),
\end{equation}
where $$\varphi_N(x) = \frac{\tN(xp) \tN(xq^{-1})
\tN(xp^{-1}q)}{\tN(xp^{-1}) \tN(xq) \tN(xpq^{-1})},$$ and $\theta_a(x)$
stands for the $\theta$--function with the multiplicative period $a$. These
relations entail similar relations for the currents $T_i(z)$.

The function $\varphi_N(x)$ can be characterized by the properties that it
is an elliptic function, which has three zeroes $u_1,u_2,u_3$, three poles
$-u_1,-u_2,-u_3$, and one of the poles is equal to $1/N$ of the
period. These properties imply that the function $\varphi_N(x)$ satisfies
the functional equation
$$\varphi_N(x) \varphi_N(xp) \ldots \varphi_N(xp^{N-1}) = 1.$$

We also show that the screening currents $S^+_i(z), i=1,\ldots,N-1$
satisfy, in the analytic continuation sense, the following relations:
\begin{equation}    \label{elrel}
S^+_i(z) S^+_j(w) = (-1)^{A_{ij}-1} \left( \frac{w}{z}
\right)^{A_{ij}-A_{ij}\beta-1} \frac{\theta_q\left(\dfrac{w}{z}
p^{A_{ij}/2}\right)}{\theta_q\left(\dfrac{z}{w}p^{A_{ij}/2}\right)}
S^+_j(w) S^+_i(z),
\end{equation}
where $(A_{ij})$ is the Cartan matrix. The screening currents $S^-_i(z),
i=1,\ldots,N-1$, satisfy the same relations with $q$ replaced by $p/q$ and
$\beta$ relaced by $1/\beta$. Moreover, we show that the screening currents
involved in the Wakimoto realization of $U_q(\sun)$ \cite{AOS}
also obey similar relations.

These elliptic relations define algebras, which are closely related to the
elliptic algebras introduced by A.~Odesskii and the first author in
\cite{FO}. Such an algebra $U_{q,p}(\N)$ can be viewed as an elliptic
deformation of the quantized universal enveloping algebra $U_q(\N)$ (where
$\N$ is the loop algebra of the nilpotent subalgebra $\n$ of $\g$),
introduced by V.~Drinfeld \cite{D}. According to \cite{FO}, the elliptic
relations of the type \eqref{elrel} imply that the screening currents
satisfy certain elliptic analogues of the quantum Serre relations from
$U_q(\N)$. We will study these relations in more detail in the next
paper. We recall that the ordinary screening charges satisfy the ordinary
quantum Serre relations \cite{BMP}, see also \cite{FF:laws}.

In this work we concentrate on the $\W$--algebras associated to $\sw_N$. In
\cite{FR} it was shown how to construct the Poisson algebra $\W_q(\g)$ and
its free field realization for general simple Lie algebra $\g$. We expect
that our results on the quantization of $\W_q(\sw_N)$ can be similarly
generalized. At the end of the paper we define the Heisenberg algebra
$\Hh_{q,p}(\g)$ and the screening currents corresponding to the general
simply-laced simple Lie algebra $\g$. We then define the algebra
$\W_{q,p}(\g)$ as the commutant of the screening charges in
$\Hh_{q,p}(\g)$. We hope that the homological methods that we used in the
study of the ordinary $\W$--algebras \cite{FF:laws} can be applied to these
quantum $\W$--algebras.

The ordinary $\W$--algebras can be obtained by the quantum Drinfeld-Sokolov
reduction from the affine algebras. We expect that the quantum
$\W$--algebras can be obtained by an analogous reduction from the quantum
affine algebras.

\subsection{} The paper is organized as follows. In Sect.~2 we recall the
results of \cite{FR} on the Poisson algebras $\W_q(\sw_N)$. In Sect.~3 we
recall the results of \cite{SKAO} on the algebra $\W_{q,p}(\sw_2)$. We
define the algebras $\W_{q,p}(\sw_N)$ in Sect.~4, and their screening
currents in Sect.~5. In Sect.~6 we derive relations in the algebra
$\W_{q,p}(\sw_N)$. In Sect.~7 we present these relations in elliptic
form. Finally, in Sect.~8 we derive the elliptic relations obeyed by the
screening currents of $\W_{q,p}(\g)$ and $U_q(\G)$.

\section{Poisson algebras $\W_q(\sw_N)$.}    \label{fr}
In this section we recall results of \cite{FR}. Let us first introduce the
Heisenberg-Poisson algebra $\Hh_q(\sw_N)$. It has generators $a_i[n],
i=1,\ldots,N-1; n\in\Z$, and relations
\begin{equation}    \label{pba}
\{ a_i[n],a_j[m] \} = h (q^{nA_{ij}/2}-q^{-nA_{ij}/2}) \delta_{n,-m},
\end{equation}
where $(A_{ij})$ is the Cartan matrix of $\sw_N$.

\begin{rem} The parameter $q$ that we use in this paper corresponds to
$q^2$ in \cite{FR}. The algebra $\W_q(\g)$ corresponds to $\W_{h/2}(\g)$ in
the notation of \cite{FR}.\qed
\end{rem}

Define now new generators $\la_i[n], i=1,\ldots,N; n\in\Z$, according to
the formula $$\la_i[n] - \la_{i+1}[n] = q^{ni/2} a_i[n], \quad \quad
i=1,\ldots,N-1; n\in\Z,$$
\begin{equation}    \label{linrel}
\sum_{i=1}^N q^{(1-i)n} \la_i[n] = 0.
\end{equation}
{}From these formulas we derive the Poisson brackets (see \cite{FR})
\begin{align}    \label{pbl1}
\{ \la_i[n],\la_i[m] \} &= - h  \frac{(1-q^n)(1-q^{n(N-1)})}{1-q^{nN}}
\delta_{n,-m}, \\    \label{pbl2}
\{ \la_i[n],\la_j[m] \} &= h \frac{(1-q^n)^2}{1-q^{nN}} q^{-n}
\delta_{n,-m}, \quad \quad i<j.
\end{align}

Introduce the generating functions
\begin{equation}    \label{lai}
\La_i(z) = q^{i-(N+1)/2} \exp \left( - \sum_{m\in\Z} \la_i[m] z^{-m}
\right).
\end{equation}
{}From \eqref{pbl1} and \eqref{pbl2}
we find:
\begin{align}    \label{pbg1}
\{ \La_i(z),\La_i(w) \} & = - h \left\{ \sum_{m\in\Z} \left( \frac{w}{z}
\right)^m \frac{(1-q^m)(1-q^{m(N-1)})}{1-q^{mN}} \right\} \La_i(z)
\La_i(w),
\\    \label{pbg2}
\{ \La_i(z),\La_j(z) \} & = h \left\{ \sum_{m\in\Z} \left(
\frac{w}{zq} \right)^m \frac{(1-q^m)^2}{1-q^{mN}} \right\} \La_i(z)
\La_j(z),
\end{align}
if $i<j$.

Now let us define generating functions $\s_i(z), i=0,\ldots,N$, whose
coefficients lie in $\HH$: $\s_0(z) = 1$, and
\begin{equation}    \label{qmap}
\s_i(z) = \sum_{1\leq j_1 < \ldots < j_i\leq N} \La_{j_1}(z) \La_{j_2}(zq)
\ldots \La_{j_{i-1}}(zq^{i-2}) \La_{j_i}(zq^{i-1}),
\end{equation}
$i=1,\ldots,N$.  Formula \eqref{linrel} implies that
$$\s_N(z) = \La_1(z) \La_2(zq) \ldots \La_N(zq^{N-1}) = 1.$$

Formula \eqref{qmap} can be rewritten succinctly as follows:
\begin{multline}    \label{sq1}
\Dq^N - \s_1(z) \Dq^{N-1} + \s_2(z) \Dq^{N-2} - \ldots + (-1)^{N-1}
\s_{N-1}(z) \Dq + (-1)^N =
\\    \label{sq2}
(\Dq - \La_1(z))(\Dq - \La_2(zq)) \ldots (\Dq - \La_N(zq^{N-1})),
\end{multline}
where $D_a$ stands for the $a$--difference operator: $$(D_a \cdot f)(x) =
f(xa).$$

In the limit $q \arr 1$ we have: $\La_i(z) = 1 - h \chi_i(z) + o(h)$ and
$\Dq = 1 - h \pa_z + o(h)$, where $h = \log q$. Hence the right hand side
of \eqref{sq1} becomes in this limit $$(-1)^N h^N (\pa_z - \chi_1(z))(\pa_z
- \chi_2(z)) \ldots (\pa_z - \chi_N(z)) + o(h^N),$$ and we obtain the
standard Miura transformation corresponding to the classical $\W$--algebra
$\W(\sw_N)$, see e.g. \cite{DS}. This shows that the generators of
$\W(\sw_N)$ can be recovered as certain linear combinations of
$\s_0(z),\ldots,\s_{N-1}(z)$ and their derivatives in the limit $q \arr 1$.

The coefficients of the series $\s_i(z), i=1,\ldots,N-1$, generate a
Poisson subalgebra $\W_q(\sw_N)$ of $\HH$. The relations between them are
as follows (see \cite{FR}):
\begin{align*}
\{ \s_i(z),\s_j(w) \} = & - h \left\{ \sum_{m\in\Z} \left(
\frac{wq^{j-i}}{z} \right)^m \frac{(1-q^{im})(1-q^{m(N-j)})}{1-q^{mN}}
\right\} \s_i(z) \s_j(w) \\ &+ h \sum_{r=1}^i \delta\left( \frac{w}{zq^r}
\right) \s_{i-r}(w) \s_{j+r}(z) \\ &- h \sum_{r=1}^i \delta\left(
\frac{wq^{j-i+r}}{z} \right) \s_{i-r}(z) \s_{j+r}(w),
\end{align*}
if $i\leq j$ and $i+j\leq N$; and
\begin{align*}
\{ \s_i(z),\s_j(w) \} = &-h \left\{ \sum_{m\in\Z} \left( \frac{wq^{j-i}}{z}
\right) \frac{(1-q^{im})(1-q^{m(N-j)})}{1-q^{mN}} \right\} \s_i(z) \s_j(w)
\\ &+ h \sum_{r=1}^{N-j} \delta\left( \frac{w}{zq^r} \right) \s_{i-r}(w)
\s_{j+r}(z) \\ &- h \sum_{r=1}^{N-j} \delta\left( \frac{wq^{j-i+r}}{z}
\right) \s_{i-r}(z) \s_{j+r}(w),
\end{align*}
if $i\leq j$ and $i+j>N$.

\begin{rem}
It is natural to define the Poisson algebra $\W_q(\sw_\infty)$ with
generators $\s_i(z)$, $i\geq 1$, and relations
\begin{align*}
\{ \s_i(z),\s_j(w) \} &= h \sum_{r=0}^i \delta\left( \frac{w}{zq^r}
\right) \s_{i-r}(w) \s_{j+r}(z) \\ &- h \sum_{r=0}^i \delta\left(
\frac{wq^{j-i+r}}{z} \right) \s_{i-r}(z) \s_{j+r}(w), \quad i\leq j.
\end{align*}\qed
\end{rem}

\section{The algebra $\W_{q,p}(\sw_2)$.}    \label{sl2}
In this section we recall the results of Shiraishi-Kubo-Awata-Odake
\cite{SKAO} on the quantum deformation of $\W_q(\sw_2)$. However, some
of our notation will be different from theirs.

Let $h,\beta$ be two complex numbers, such that neither $h$ nor $h\beta$
belongs to $2\pi i {\Bbb Q}$. Set $q=e^h$ and $p=e^{h(1-\beta)}$. We will
use this notation throughout the paper.

Let $\Hh'_{q,p}(\sw_2)$ be the Heisenberg algebra with generators $\la[n],
n\in\Z$, and relations:
\begin{equation}    \label{hprel}
\left[ \la[n],\la[m] \right] = -\frac{1}{n}
\frac{(1-q^n)(1-(p/q)^n)}{1+p^n} \delta_{n,-m}.
\end{equation}
In the limit $\beta \arr 0$, in which $p \arr q$ we can recover the Poisson
bracket \eqref{hrel} as the $\beta$--linear term of the bracket
\eqref{hprel}.

For $\mu \in \C$, let $\pi_\mu$ be the Fock representation of the algebra
$\Hh'_{q,p}(\sw_2)$, which is generated by a vector $v_\mu$, such that
$\la[n] v_\mu = 0, n>0$, and $\la[0] v_\mu = \mu v_\mu$.

Let $$\Hh_{q,p}(\sw_2) = \underset{\longleftarrow}{\lim} \;
\Hh'_{q,p}(\sw_2)/I_n, \quad \quad n>0,$$ where $I_n$ is the left ideal of
$\Hh'_{q,p}(\sw_2)$ generated by all polynomials in $\la[m], m>0$, of
degrees greater than or equal to $n$ ($\deg \la[m] = m$). By definition,
the action of $\Hh_{q,p}(\sw_2)$ on the modules $\pi_\mu$ is well-defined.

Introduce the generating function $$\La(z) = p^{-1/2} q^{-\la[0]} :\exp
\left( - \sum_{m\neq 0} \la[m] z^{-m} \right):,$$ where columns stand for
the standard normal ordering. Now define the power series $T(z) =
\sum_{m\in\Z} T[m] z^{-m}$ by the formula
\begin{equation}    \label{qvirrel}
T(z) = :\La(zp^{1/2}): + :\La(zp^{-1/2})^{-1}:.
\end{equation}
The coefficients $T[n]$ of the power series $T(z)$ belong to
$\Hh_{q,p}(\sw_2)$. They satisfy the following relations \cite{SKAO}:
\begin{multline}    \label{skao}
\sum_{l=0}^\infty f_l \left(T[n-l] T[m+l] - T[m-l] T[n+l]\right) \\ =
\frac{(1-q)(1-p/q)}{1-p} (p^{-n}-p^n) \delta_{n,-m},
\end{multline}
where $f_l$'s are given by the generating function
\begin{equation}    \label{fz}
f(x) = \sum_{l=0}^\infty f_l x^l = \exp \left( \sum_{m=1}^\infty
\frac{1}{m} \frac{(1-q^m)(1-(p/q)^m)}{1+p^m} x^m \right).
\end{equation}

In the limit $q \arr 1$, formulas \eqref{qvirrel} and \eqref{skao} become
formulas \eqref{ff} and \eqref{virrel}, respectively.

Introduce an additional operator $Q$, such that $\left[ \la[n],Q \right] =
\beta \delta_{n,0}$. The operator $e^{\al Q}, \al\in\C$, acts from $\pi_\mu$
to $\pi_{\mu+\al\beta}$ by sending $v_\mu$ to $v_{\mu+\al\beta}$. In
\cite{SKAO} two screening currents were constructed:
\begin{align}    \label{splus}
S^+(z) &= e^{Q} z^{s^+[0]} :\exp \left( \sum_{m\neq 0} s^+[m] z^{-m}
\right):,\\ \label{sminus} S^-(z) &= e^{-Q/\beta} z^{-s^-[0]} :\exp \left(
- \sum_{m\neq 0} s^-[m] z^{-m} \right):,
\end{align}
where $$s^+[m] = \frac{1+p^m}{q^{-m}-1} \la[m], \quad m \neq 0, s^+[0] = 2
\la[0],$$ $$s^-[m] = \frac{1+p^m}{(q/p)^m-1} \la[m], \quad m \neq 0, \quad
\quad s^-[0] = 2 \la[0]/\beta.$$ The Fourier coefficients of $S^+(z)$ act
from $\pi_\mu$ to $\pi_{\mu+\beta}$, and the Fourier coefficients of
$S^-(z)$ act from $\pi_\mu$ to $\pi_{\mu-1}$.

They satisfy \cite{SKAO}:
$$[T[n],S^+(w)] = {\cal D}_q C^+_n(w), \quad \quad [T[n],S^-(w)] = {\cal
D}_{p/q} C^-_n(w),$$ where $C^\pm_n(w)$ are certain operator--valued power
series, and $$[{\cal D}_a \cdot f](x) = \frac{f(x) - f(xa)}{x(1-a)}.$$ This
implies that $T[n], n\in\Z$, commute with the screening charges $\int
S^\pm(z) dz$, whenever they are well-defined \cite{SKAO}. In the limit $q
\arr 1$, those become the two screening charges of the Virasoro algebra.

\section{The algebra $\W_{q,p}(\sw_N)$.}
\subsection{Heisenberg algebra.}    \label{heisenberg}
Let $\Hh'_{q,p}(\sw_N)$ be the Heisenberg algebra with generators $a_i[n],
i=1,\ldots,N-1; n\in\Z$, and relations
\begin{equation}    \label{ai}
\left[ a_i[n],a_j[m] \right] = \frac{1}{n}
\frac{(1-q^n)(p^{A_{ij}n/2}-p^{-A_{ij}n/2})(1-(p/q)^n)}{1-p^n}
\delta_{n,-m},
\end{equation}
This formula was guessed from the commutation relations \eqref{hprel} in
the case of $\sw_2$, which follow from \cite{SKAO}, and from the condition
that in the limit $\beta \arr 0$ the $\beta$--linear term should give us
the Poisson bracket \eqref{pba}.

For each weight $\mu$ of the Cartan subalgebra of $\sw_N$, let $\pi_\mu$ be
the Fock representation of $\Hh'_{q,p}(\sw_N)$ generated by a vector
$v_\mu$, such that $a_i[n] v_\mu = 0, n>0$, and $a_i[0] v_\mu =
\mu(\al_i^\vee) v_\mu$, where $\al_i^\vee$ is the $i$th coroot of $\sw_N$.

Let $\Hh_{q,p}(\sw_N)$ be the completion of $\Hh'(\sw_N)$ defined in the
same way as in the case of $\sw_2$, see \secref{sl2}. The algebra
$\Hh_{q,p}(\sw_N)$ acts on the modules $\pi_\mu$.

Introduce new generators $\la_i[n]$ of $\Hh_{q,p}(\sw_N)$ by the formulas
\begin{equation}    \label{ladef}
\la_i[n] - \la_{i+1}[n] = p^{ni/2} a_i[n], \quad \quad i=1,\ldots,N-1;
n\in\Z,
\end{equation}
\begin{equation}    \label{linrelp}
\sum_{i=1}^N p^{(1-i)n} \la_i[n] = 0.
\end{equation}
{}From these formulas and \eqref{ai} we derive the commutation relations
between them:
\begin{align}    \label{bl1}
\left[ \la_i[n],\la_i[m] \right] &= - \frac{1}{n}
\frac{(1-q^n)(1-p^{n(N-1)})(1-(p/q)^n)}{1-p^{nN}}
\delta_{n,-m}, \\    \label{bl2}
\left[ \la_i[n],\la_j[m] \right] &= \frac{1}{n}
\frac{(1-q^n)(1-p^n)(1-(p/q)^n)}{1-p^{nN}} p^{-n} \delta_{n,-m}, \quad \quad
i<j.
\end{align}

Let us introduce the power series
\begin{equation}    \label{li}
\Lambda_i(z) = p^{i-(N+1)/2} q^{-\la_i[0]} :\exp \left( - \sum_{m\neq 0}
\la_i[m] z^{-m} \right):.
\end{equation}

We can compute the operator product expansions (OPEs) of these power series
using the following lemma. Introduce the notation
$$(x|\al_1,\ldots,\al_k;t)_\infty = \prod_{i=1}^k \prod_{n=0}^\infty
(1-\al_i t^n).$$

\begin{lem}    \label{ope}
Let $b[n], n\in\Z$, and $c[n], n\in\Z$, satisfy commutation
relations:
$$\left[b[n],c[m]\right] = - \frac{1}{n(1-t)} \left( \sum_{i=1}^k \al_i^n -
\sum_{j=1}^l \beta_j^n \right) \delta_{n,-m},$$ where $|t|<1$. Then for
$|z| > \underset{i,j}{\on{max}} \left\{ |\al_i|,|\beta_j|\right\} |w|$ the
composition
$$:\exp \left( \sum_{n\in\Z} b[n] z^{-n} \right): :\exp \left(
\sum_{n\in\Z} c[n] z^{-n} \right):$$ acting on each module $\pi_\mu$ exists
and is equal to $$\frac{\ds \left(\dfrac{w}{z}| \al_1,\ldots,\al_k
\right)_\infty}{\ds \left(\dfrac{w}{z}| \beta_1,\ldots, \beta_l
\right)_\infty} :\exp \left( \sum_{n\in\Z} b[n] z^{-n} \right) \exp \left(
\sum_{n\in\Z} c[n] z^{-n} \right):.$$
\end{lem}

\noindent{\em Proof.} Direct computation based on formula $$\exp \left( -
\sum_{n>0} \frac{x^n}{n} \right) = 1-x.$$\qed

Let us assume that $|p|<1$ and $|z| \gg |w|$; more precisely, it suffices
that $|z|>|w|pq^{-1}$, and $|z|>|z|q$. Then we find from formula
\eqref{bl1} and \lemref{ope}:
\begin{equation}    \label{opel1}
\La_i(z) \La_i(w) = \frac{\left(\dfrac{w}{z}|
1,p,p^{N-1}q,p^Nq^{-1};p^N\right)_\infty}{\left(\dfrac{w}{z}|
q,pq^{-1},p^{N-1},p^N;p^N\right)_\infty} :\La_i(z) \La_i(w):.
\end{equation}

In the same way we obtain:
\begin{align}    \label{opel2}
\La_i(z) \La_j(w) & = \frac{\left(\dfrac{w}{z}|
qp^{-1},q^{-1},p;p^N\right)_\infty}{\left(\dfrac{w}{z}|
p^{-1},q,pq^{-1};p^N\right)_\infty} :\La_i(z)
\La_j(w):, \quad i<j,\\    \label{opel3}
\La_i(z) \La_j(w) & = \frac{\left(\dfrac{w}{z}|
p^{N-1}q,p^Nq^{-1},p^{N+1};p^N\right)_\infty}{\left(\dfrac{w}{z}|
p^{N-1},p^Nq,p^{N+1}q^{-1};p^N\right)_\infty} :\La_i(z) \La_j(w):, \quad
i>j,
\end{align}
when $|z|\gg|w|$.

\begin{rem}    \label{power}
When $|p|<1$, the functions appearing in the right hand side of formulas
\eqref{opel1}--\eqref{opel3} are power series in $w/z$, whose coefficients
are rational functions in $p$.\qed
\end{rem}

\subsection{Definition of the quantum $\W$--algebra.}
Now we define generating functions $T_i(z), i=0,\ldots,N$, whose
coefficients lie in $\HH$: $T_0(z) = 1$, and
\begin{equation}    \label{pmap}
T_i(z) = \sum_{1\leq j_1 < \ldots < j_i\leq N} :\La_{j_1}(z) \La_{j_2}(zp)
\ldots \La_{j_{i-1}}(zp^{i-2}) \La_{j_i}(zp^{i-1}):,
\end{equation}
$i=1,\ldots,N$. Formula \eqref{linrelp} implies that
$$T_N(z) = :\La_1(z) \La_2(zp) \ldots \La_N(zp^{N-1}): = 1.$$

Formula \eqref{pmap} can be rewritten as follows:
\begin{multline}    \label{sp1}
\Dp^N - T_1(z) \Dp^{N-1} + T_2(z) \Dp^{N-2} -
\ldots + (-1)^{N-1} T_{N-1}(z) \Dp + (-1)^N = \\
:(\Dp - \La_1(z))(\Dp - \La_2(zp)) \ldots (\Dp - \La_N(zp^{N-1})):.
\end{multline}
Formulas \eqref{pmap} and \eqref{sp1} are quantum deformations of formulas
\eqref{qmap} and \eqref{sq1}--\eqref{sq2}.

We define the algebra $\W_{q,p}(\sw_N)$ as the subalgebra of $\HH$
generated by the Fourier coefficients of the power series $T_i(z),
i=1,\ldots,N-1$, given by formula \eqref{pmap}. It is clear from the
definition that in the limit $\beta \arr 0$, i.e. $p \arr q$, the
algebra $\W_{q,p}(\sw_N)$ becomes the Poisson algebra $\W_q(\sw_N)$ defined
in \cite{FR}, see \secref{fr}.

\begin{rem} The currents $\Lambda(z)$ and $T(z)$ that were used in the
definition of $\W_{q,p}(\sw_2)$ in \secref{sl2} correspond to
$\Lambda_1(z)$ and $T_1(zp^{1/2})$, respectively.\qed
\end{rem}

Let us fix $\beta$ and consider the limit $q \arr 1$ with
$p=q^{1-\beta}$. Then we have: $\La_i(z) = 1 - h \chi_i(z) + o(h)$ and $\Dp
= 1 - h(1-\beta) \pa_z + o(h)$, where $h = \log q$. Hence the right hand
side of \eqref{sp1} becomes in this limit $$(-1)^N h^N \left(
(1-\beta)\pa_z - \chi_1(z))((1-\beta)\pa_z - \chi_2(z)) \ldots
((1-\beta)\pa_z - \chi_N(z) \right) + o(h^N),$$ and we obtain the normally
ordered Miura transformation corresponding to the $\W$--algebra of $\sw_N$,
introduced by Fateev and Lukyanov \cite{FL}. In the notation of
\cite{FF:laws}, this algebra is $\W_{\sqrt{\beta}}(\sw_N)$ with central
charge $(N-1) - N(N+1)(1-\beta)^2/\beta$. Thus, in the limit $q \arr 1$,
the algebra $\W_{q,p}(\sw_N)$ becomes isomorphic to
$\W_{\sqrt{\beta}}(\sw_N)$. The generating currents of
$\W_{\sqrt{\beta}}(\sw_N)$ can be recovered as certain linear combinations
of $T_0(z),\ldots,T_{N-1}(z)$ and their derivatives in the limit $q \arr
1$.

\section{Screening currents for $\W_{q,p}(\sw_N)$.}
Introduce operators $Q_i, i=1,\ldots,N-1$, which satisfy commutation
relations $[a_i[n],Q_j] = A_{ij} \beta \delta_{n,0}$. The operators
$e^{Q_i}$ act from $\pi_\mu$ to $\pi_{\mu+\beta\al_i}$.

Now we can define the screening currents as the generating functions
\begin{align}    \label{s+}
S_i^+(z) & = e^{Q_i} z^{s^+_i[0]} :\exp \left( \sum_{m\neq 0} s^+_i(m)
z^{-m} \right):,\\    \label{s-}
S_i^-(z) & = e^{-Q_i/\beta} z^{-s^-_i[0]} :\exp \left( - \sum_{m\neq 0}
s^-_i(m) z^{-m} \right):,
\end{align}
where
\begin{equation}    \label{sm+}
s^+_i[m] = \frac{a_i[m]}{q^{-m}-1}, \quad m\neq 0, \quad \quad
s^+_i[0] = a_i[0],
\end{equation}
\begin{equation}    \label{sm-}
s^-_i[m] = - \frac{a_i[m]}{(q/p)^m-1}, \quad m\neq 0,
\quad \quad s^-_i[0] = a_i[0]/\beta
\end{equation}
(compare with \eqref{splus} and \eqref{sminus}).

Let
\begin{equation}    \label{aseries}
A_i(z) = q^{-a[0]} :\exp \left( - \sum_{m\neq 0} a_i[m] z^{-m}
\right):.
\end{equation}
Then we have:
\begin{equation}    \label{skr1}
A_i(z) = :S^+_i(z) S^+_i(zq)^{-1}:,
\end{equation}
and
\begin{equation}    \label{skr2}
A_i(z) = :S^-_i(z) S^-_i(zp/q)^{-1}:.
\end{equation}

Formulas \eqref{skr1} and \eqref{skr2} show that the screening currents are
solutions of the following difference equations:
\begin{align*}
D_q S^+_i(z) & = :A_i(z)^{-1} S^+_i(z):,\\ D_{p/q} S^-_i(z) & =
:A_i(z)^{-1} S^-_i(z):.
\end{align*}
In the limit $q \arr 1$ they become the differential equations defining the
ordinary screening currents.

We also have from \eqref{ladef}:
\begin{equation}    \label{ala}
A_i(z) = p:\La_i(zp^{i/2}) \La_{i+1}(zp^{i/2})^{-1}:.
\end{equation}

\begin{thm} {\em The screening currents commute with the algebra
$\W_{q,p}(\sw_N)$ up to a total difference. More precisely, for any $A \in
\W_{q,p}(\sw_N)$ we have: $$[A,S^+_i(w)] = {\cal D}_q C^+_i(w), \quad \quad
[A,S^-_i(w)] = {\cal D}_{p/q} C^-_i(w),$$ where $C^\pm(w)$ are certain
operator-valued power series, and} $$[{\cal D}_a f](x) = \frac{f(x) -
f(xa)}{x(1-a)}.$$
\end{thm}

\noindent{\em Proof.} Let us consider the case of the screening currents
$S^+_i(z)$; the case of $S^-_i(z)$ can be treated in the same way.

Consider the difference operator \eqref{sp1}. We want to prove that each
term of this operator has the property that all of its Fourier coefficients
commutes with $S^+_i(z)$ up to a total ${\cal D}_q$--difference.

{}From formulas \eqref{sm+}, \eqref{ladef}, \eqref{bl1}, \eqref{bl2} we
obtain the following commutation relations:
$$\left[ \la_i[n],s^+_i[m] \right] = \frac{1}{n} p^{n(i/2-1)} (1-(p/q)^n)
\delta_{n,-m},$$
$$\left[ \la_{i+1}[n],s^+_i[m] \right] = - \frac{1}{n} p^{ni/2} (1-(p/q)^n)
\delta_{n,-m},$$
$$\left[ \la_j[n],s^+_i[m] \right] = 0, \quad \quad j\neq i,i+1.$$

{}From these commutation relations we derive the following OPEs
(cf. \lemref{ope}):
$$\La_i(z) S^+_i(w) =
\frac{p(z-wp^{i/2-1})}{q(z-wp^{i/2}q^{-1})} :\La_i(z)
S^+_i(w):, \quad \quad |z|\gg|w|,$$ $$S^+_i(w) \La_i(z) =
\frac{p(z-wp^{i/2-1})}{q(z-wp^{i/2}q^{-1})} :\La_i(z)
S^+_i(w):, \quad \quad |w|\gg|z|,$$ $$\La_{i+1}(z) S^+_i(w) =
\frac{q(z-wp^{i/2+1}q^{-1})}{p(z-wp^{i/2})} :\La_{i+1}(z)
S^+_i(w):, \quad \quad |z|\gg|w|,$$ $$S^+_i(w) \La_{i+1}(z) =
\frac{q(z-wp^{i/2+1}q^{-1})}{p(z-wp^{i/2})} :\La_{i+1}(z)
S^+_i(w):, \quad \quad |w|\gg|z|,$$ and $$\La_j(z) S^+_i(w) = :\La_j(z)
S^+_i(w):, \quad \quad \forall z,w,$$ if $j\neq i,i+1$.

The last formula means that $S^+_i(w)$ commutes with all Fourier
coefficients of $\La_j(z)$ if $j\neq i,i+1$. Therefore it is sufficient to
consider the OPE between the factor
\begin{multline*}
:(\Dp - \La_i(zp^{i-1}))(\Dp - \La_{i+1}(zp^i)): \\ = \Dp^2 -
\left(\La_i(zp^{i-1}) + \La_{i+1}(zp^{i-1})\right) \Dp + :\La_i(zp^{i-1})
\La_{i+1}(zp^i):
\end{multline*}
in formula \eqref{sp1} and $S^+_i(w)$. We have to show that all Fourier
coefficients of each of the terms commute with $S^+_i(w)$ up to a total
difference.

For the term $:\La_i(zp^{i-1}) \La_{i+1}(zp^i):$ we have according to the
OPEs above:
$$:\La_i(zp^{i-1}) \La_{i+1}(zp^i): S^+_i(w) = :\La_i(zp^{i-1})
\La_{i+1}(zp^i) S^+_i(w):,$$ which means that all Fourier coefficients of
$:\La_i(zp^{i-1}) \La_{i+1}(zp^i):$ commute with $S^+_i(w)$.

Now consider the linear term $\La_i(zp^{i-1}) + \La_{i+1}(zp^{i-1})$. We
have according to the OPEs above:
$$\left(\La_i(zp^{i-1}) + \La_{i+1}(zp^{i-1})\right) S^+_i(w) =$$
\begin{multline}    \label{qdiff}
\frac{p(z-wp^{-i/2})}{q(z-wp^{-i/2+1}q^{-1})} :\La_i(zp^{i-1})
S^+_i(w): \\ + \frac{q(z-wp^{i/2+2}q^{-1})}{p(z-wp^{-i/2+1})}
:\La_{i+1}(zp^{i-1}) S^+_i(w):,
\end{multline}
for $|z|\gg|w|$, and the same formula for the product in the opposite order
for $|w|\gg|z|$. Therefore we can compute the commutator
$$\left[ \int \left(\La_i(zp^{i-1}) + \La_{i+1}(zp^{i-1})\right) z^n
dz,S^+_i(w) \right]$$ by evaluating the residues in the right hand side of
\eqref{qdiff}. We find that this commutator is equal to
\begin{multline*}
\left( p/q - 1 \right) \left[ :\La_i(wp^{i/2}q^{-1}) S^+_i(w):
(wp^{-i/2+1}q^{-1})^{n+1} \right. \\ \left. - p^{-1}q :\La_{i+1}(wp^{i/2})
S^+_i(w): (wp^{-i/2+1})^{n+1} \right]
\end{multline*}
$$= {\cal D}_q \left[ p^{i/2-1} \left( 1-q \right) \left( p - q \right)
:\La_i(wp^{i/2}q^{-1}) S^+_i(w): (wp^{-i/2+1}q^{-1})^{n+2} \right],$$
because by formulas \eqref{skr1} and \eqref{ala} $$:\La_i(wp^{i/2})
S^+_i(wq): = p^{-1}:\La_{i+1}(wp^{i/2}) S^+_i(w):.$$ This completes the
proof.\qed

\begin{cor}    \label{residue}
Any element of the algebra $\W_{q,p}(\sw_N)$ commutes with the operators
$\int S^\pm_i(z) dz$ acting on the Fock representations, whenever they are
well-defined.
\end{cor}

In the limit $q \arr 1$, the operators $\int S^\pm_i(z) dz$ become the
ordinary screening charges $\int e^{\al_\pm \phi_i(z)} dz$, where $\al_\pm
= \pm \beta^{\pm 1/2}$.

\corref{residue} implies that one can construct intertwining operators
between the Fock representations $\pi_\mu$ of $\W_{q,p}(\sw_N)$, and hence
singular vectors, by integrating products of the screening currents over
suitable cycles. For the ordinary $\W$--algebra of $\sw_N$, the screening
charges satisfy quantum Serre relations, and the integration cycles
correspond to the singular vectors in the Verma modules over the quantum
group $U_q(\sw_N)$, see \cite{BMP,SV,V,FF:laws}. We expect an analogous
structure for the $\W_{q,p}(\sw_N)$ screening charges. We will return to
this question in our next paper.

The singular vectors in the Fock representations of $\W_{q,p}(\sw_N)$
should coincide with the Macdonald symmetric functions corresponding to the
Young diagrams with $N-1$ rectangles, as was pointed out in \cite{SKAO}. In
the case of $\sw_2$ formulas for these singular vectors were given in
\cite{SKAO,AOS:bos}.

\section{Relations in $\W_{q,p}(\sw_N)$.}
\subsection{Relations between $T_1(z)$ and $T_m(w)$.}

Let us again assume that $|p|<1$. Introduce the formal power series
$f_{m,N}(x)$ by the formula
\begin{equation}    \label{fmn}
f_{m,N}(x) = \dfrac{\left(x|
p^{m-1}q,p^mq^{-1},p^{N-1},p^N;p^N\right)_\infty}{\left(x|
p^{m-1},p^m,p^{N-1}q,p^Nq^{-1};p^N\right)_\infty}.
\end{equation}

The function $f_{m,N}(x)$ is a very-well-poised basic hypergeometric series
$$_6 \phi_5 \left[ \begin{array}{llllll} x, & x^{1/2} p^N, & - x^{1/2} p^N,
& p^{N-m}, & pq^{-1}, & q \\ {} & x^{1/2}, & - x^{1/2} & xp^m, & xp^{N-1}q,
& xp^Nq^{-1} \end{array} ; p^N, x \right],$$ see formula (2.7.1) of
\cite{GR}.

\begin{thm}
{\em The formal power series $T_1(z)$ and $T_m(w)$ satisfy the following
relations:}
\begin{multline}    \label{m1}
f_{m,N} \left(\frac{w}{z}\right) T_1(z) T_m(w) -
f_{m,N} \left(\frac{z}{w}\right) T_m(w) T_1(z) \\
= \frac{(1-q)(1-p/q)}{1-p} \left( \delta \left( \dfrac{w}{zp} \right)
T_{m+1}(z) - \delta \left( \dfrac{wp^m}{z} \right) T_{m+1}(w) \right).
\end{multline}
\end{thm}

\noindent{\em Proof.}  Using the OPEs \eqref{opel1}--\eqref{opel3}, we
obtain that when $|z|\gg|w|$ $$\La_i(z) :\La_{j_1}(w) \La_{j_2}(wp)
\ldots \La_{j_m}(wp^{m-1}):$$ is equal to
$$f_{m,N}\left(\frac{w}{z}\right)^{-1} :\La_i(z) \La_{j_1}(w) \La_{j_2}(wp)
\ldots \La_{j_m}(wp^{m-1}):,$$ if $i=j_k$ for some $k \in \{ 1,\ldots,m
\}$; and
$$f_{m,N}\left(\frac{w}{z}\right)^{-1} \frac{(z-w p^{k-1}q)
(z-wp^kq^{-1})}{(z-wp^{k-1})(z-wp^k)} :\La_i(z) \La_{j_1}(w) \La_{j_2}(wp)
\ldots \La_{j_m}(wp^{m-1}):,$$ if $j_k<i<j_{k+1}$. Here and below the case
$i<j_1$ corresponds to $k=0$ and the case $i>j_m$ corresponds to $k=m$.

On the other hand, when $|w|\gg|z|$, $$:\La_{j_1}(w)
\La_{j_2}(wp) \ldots \La_{j_m}(wp^{m-1}): \La_i(z)$$ is equal to
$$f_{m,N}\left(\frac{z}{w}\right)^{-1} :\La_{j_1}(w) \La_{j_2}(wp) \ldots
\La_{j_m}(wp^{m-1}) \La_i(z):,$$ if $i=j_k$ for some $k \in \{ 1,\ldots,m
\}$; and $$f_{m,N}\left(\frac{z}{w}\right)^{-1} \frac{(z-w p^{k-1}q)
(z-wp^kq^{-1})}{(z-wp^{k-1})(z-wp^k)} :\La_{j_1}(w) \La_{j_2}(wp) \ldots
\La_{j_m}(wp^{m-1}) \La_i(z):,$$ if $j_k<i<j_{k+1}$.

Since the normally ordered product does not depend on the order of the
factors, we conclude that the analytic continuations of
$$f_{m,N}\left(\frac{w}{z}\right) \La_i(z) :\La_{j_1}(w) \ldots
\La_{j_m}(wp^{m-1}):$$ and
$$f_{m,N}\left(\frac{z}{w}\right) :\La_{j_1}(w) \ldots
\La_{j_m}(wp^{m-1}): \La_i(z)$$ coincide.

Therefore
\begin{multline*}
\int_{C_R} f_{m,N}\left(\frac{w}{z}\right) \La_i(z)
:\La_{j_1}(w) \La_{j_2}(wp) \ldots \La_{j_m}(wp^{m-1}): z^n dz \\ -
\int_{C_r} f_{m,N}\left(\frac{z}{w}\right) :\La_{j_1}(w)
\La_{j_2}(wp) \ldots \La_{j_m}(wp^{m-1}): \La_i(z) z^n dz
\end{multline*}
where $C_R$ and $C_r$ are circles on the $z$ plane of radii $R\gg|w|$ and
$r\ll|w|$, respectively, is equal to $0$ if $i=j_k$ for some $k \in \{
1,\ldots,m \}$; and the sum of the residues of $$\frac{(z-w p^{k-1}q)
(z-wp^kq^{-1})}{(z-wp^{k-1})(z-wp^k)} :\La_i(z) \La_{j_1}(w) \ldots
\La_{j_m}(wp^{m-1}):,$$ if $j_k<i<j_{k+1}$.

But the latter is equal to $(1-q)(1-p/q)/(1-p)$ times
$$\left( :\La_{j_1}(w) \ldots \La_{j_k}(wp^{k-1}) \La_i(wp^{k-1})
\La_{j_{k+1}}(wp^k) \ldots \La_{j_m}(wp^{m-1}): w^{n+1} p^{(n+1)(k-1)}
\right.$$
$$\left. - :\La_{j_1}(w) \ldots \La_{j_k}(wp^{k-1}) \La_i(wp^k)
\La_{j_{k+1}}(wp^k) \ldots \La_{j_m}(wp^{m-1}): w^{n+1} p^{(n+1)k}
\right).$$ After summation over $j_1 < j_2 < \ldots < j_m$, all of these
terms will cancel out except for $$\frac{(1-q)(1-p/q)}{1-p} :\La_i(wp^{-1})
\La_{j_1}(w) \ldots \La_{j_m}(wp^{m-1}): w^{n+l+1} p^{-(n+1)}$$ with
$i<j_1$; and $$- \frac{(1-q)(1-p/q)}{1-p} :\La_{j_1}(w) \ldots
\La_{j_m}(wp^{m-1}) \La_i(wp^m): w^{n+l+1} p^{(n+1)m}$$ with $i>j_m$. This
gives us formula \eqref{m1}.\qed

In the limit $p \arr q$ formula \eqref{m1} gives the Poisson bracket
between $\s_1(z)$ and $\s_m(w)$ from \cite{FR}, see \secref{fr}.

Formula \eqref{m1} shows that the Fourier coefficients $T_1[n]$ of the
power series $T_1(z)$ generate the algebra $\W_{q,p}(\sw_N)$. In
particular, $T_i(z)$ can be written as a degree $i$ expression in $T_1[n],
n\in\Z$.

One can also derive similar relations between $T_i(z)$ and $T_j(w)$ with
$i,j>1$. These relations are quadratic, and involve products of
$T_{i-r}(z)$ and $T_{j+r}(w)$, where $r=1,\ldots,i-1$, if $i<j$. In the
limit $p \arr q$ they give the Poisson brackets between $\s_i(z)$ and
$\s_j(w)$ from \cite{FR}, see \secref{fr}.

Let us define analogues of the Verma modules over the algebra
$\W_{q,p}(\sw_N)$; in the case of $\sw_2$ this has been done in
\cite{SKAO}.

Although the $0$th Fourier coefficients $T_i[0]$ of the series $T_i(z)$ do
not commute with each other, they commute modulo the ideal generated by
$T_i[n], i=1,\ldots,N-1; n>0$. We can therefore define a Verma module
$M_{\gamma_1,\ldots,\gamma_{N-1}}$ as a $\W_{q,p}(\sw_N)$--module generated
by a vector $v_{\gamma_1,\ldots,\gamma_{N-1}}$, such that $T_i[n]
v_{\gamma_1,\ldots,\gamma_{N-1}} = 0$, if $n>0$, and $T_i[0]
v_{\gamma_1,\ldots,\gamma_{N-1}} = \gamma_i
v_{\gamma_1,\ldots,\gamma_{N-1}}$. The relations in $\W_{q,p}(\sw_N)$ imply
that the module $M_{\gamma_1,\ldots,\gamma_{N-1}}$ has a PBW basis which
consists of lexicographically ordered monomials in $T_i[n], n<0$, applied
to the highest weight vector $v_{\gamma_1,\ldots,\gamma_{N-1}}$.

\subsection{Relations in $\W_{q,p}(\sw_2)$.}
In this case the relations are:
\begin{multline}    \label{wqpsl2}
f_{1,2} \left(\frac{w}{z}\right) T_1(z) T_1(w) - f_{1,2}
\left(\frac{z}{w}\right) T_1(w) T_1(z) \\ = \frac{(1-q)(1-p/q)}{1-p} \left(
\delta \left( \dfrac{w}{zp} \right) - \delta \left( \dfrac{wp}{z} \right)
\right).
\end{multline}
These relations are equivalent to the relations of Shiraishi,
e.a. \cite{SKAO} given by formula \eqref{skao}, because their $f(x)$
coincides with $f_{1,2}(x)$ given by \eqref{fmn}. Formula \eqref{fmn} can
be simplified in this case: $$f_{1,2}(x) = \frac{1}{1-x}
\frac{(x|q,pq^{-1};p^2)_\infty}{(x|pq,p^2q^{-1};p^2)_\infty}.$$

\subsection{Relations in $\W_{q,p}(\sw_3)$.}
In this case we have:
\begin{align*}
f_{1,3}(x) & = \frac{(x| q,pq^{-1},p^2,p^3;p^3)_\infty}{(x|
1,p,p^2q,p^3q^{-1};p^3)_\infty},\\
f_{2,3}(x) & = \frac{(x| pq,p^2q^{-1},p^3;p^3)_\infty}{(x|
p,p^2q,p^3q^{-1};p^3)_\infty}.
\end{align*}
The relations are the following:
\begin{multline*}
f_{1,3} \left(\frac{w}{z}\right) T_1(z) T_1(w) -
f_{1,3} \left(\frac{z}{w}\right) T_1(w) T_1(z) \\
= \frac{(1-q)(1-p/q)}{1-p} \left( \delta \left( \dfrac{w}{zp} \right)
T_2(z) - \delta \left( \dfrac{wp}{z} \right) T_2(w) \right);
\end{multline*}
\begin{multline*}
f_{2,3} \left(\frac{w}{z}\right) T_1(z) T_2(w) -
f_{2,3} \left(\frac{z}{w}\right) T_2(w) T_1(z) \\
= \frac{(1-q)(1-p/q)}{1-p} \left( \delta \left( \dfrac{w}{zp} \right)
 - \delta \left( \dfrac{wp^2}{z} \right) \right);
\end{multline*}
\begin{multline*}
f_{1,3} \left(\frac{w}{z}\right) T_2(z) T_2(w) -
f_{1,3} \left(\frac{z}{w}\right) T_2(w) T_2(z) \\
= \frac{(1-q)(1-p/q)}{1-p} \left( \delta \left( \dfrac{w}{zp} \right)
T_1(w) - \delta \left( \dfrac{wp}{z} \right) T_1(z) \right).
\end{multline*}
In the limit $p \arr q$ they become the relations in $\W_q(\sw_3)$
described in \cite{FR}.

\section{Relations in elliptic form.}
We recall that $q$ and $p$ are assumed to be generic with $|p|<1$.

\subsection{The case of $\sw_2$.}
Consider the OPE given by formula \eqref{opel1}:
\begin{equation}    \label{conti1}
\La(z) \La(w) = f_{1,2}\left(\frac{w}{z}\right)^{-1} :\La(z) \La(w):,
\end{equation}
where $$f_{1,2}\left(\frac{w}{z}\right) = \frac{\left(\dfrac{w}{z}|
q,pq^{-1},p^2;p^2 \right)_\infty}{\left(\dfrac{w}{z}| 1,pq,p^2 q^{-1};p^2
\right)_\infty}.$$ Formula \eqref{conti1} is valid for $|z|\gg|w|$ (see
\lemref{ope}) and it shows that the composition $\La(z) \La(w)$ can be
analytically continued to a meromorphic operator-valued function on $\C
\times \C$, given by the right hand side of the formula.

Likewise, the composition $\La(w) \La(z)$ converges when $|w|\gg|z|$, and we
have:
\begin{equation}    \label{conti2}
\La(w) \La(z) = f_{1,2}\left(\frac{z}{w}\right)^{-1} :\La(w) \La(z):.
\end{equation}

Since $:\La(z) \La(w): = :\La(w) \La(z):$, by definition of the normal
ordering, we obtain from formulas \eqref{conti1} and \eqref{conti2} the
following relation on the analytic continuations:
\begin{equation}    \label{equality}
\La(z) \La(w) = \varphi \left( \dfrac{w}{z} \right) \La(w) \La(z),
\end{equation}
where
\begin{equation}    \label{g}
\varphi(x) = \frac{f_{1,2}(x^{-1})}{f_{1,2}(x)} =
\frac{\theta_{p^2}(x)\theta_{p^2}(xpq)
\theta_{p^2}(xp^2q^{-1})}{\theta_{p^2} (xq)
\theta_{p^2}(xpq^{-1})\theta_{p^2}(xp^2)}.
\end{equation}
and
$$\theta_a(x) = \prod_{n=0}^\infty (1-xa^n)
\prod_{n=1}^\infty (1-x^{-1}a^n) \prod_{n=1}^\infty (1-a^n).$$

We can also write: $$\varphi(x) = \frac{\theta_{p^2}(xp)
\theta_{p^2}(xq^{-1}) \theta_{p^2}(xp^{-1}q)}{\theta_{p^2}(xp^{-1})
\theta_{p^2}(xq) \theta_{p^2}(xpq^{-1})}.$$

Formula \eqref{equality} can be rewritten in a more symmetric form:
\begin{equation}    \label{symform}
\gamma\left(\frac{w}{z}\right) \La(z) \La(w) = \gamma\left(\frac{z}{w}\right)
\La(w) \La(z),
\end{equation}
in the sense of analytic continuation, where
\begin{equation}    \label{h}
\gamma(x) = \frac{\theta_{p^2}(x) \theta_{p^2}(xq)}{\theta_{p^2}(xp)
\theta_{p^2}(xp^{-1}q)}.
\end{equation}

Note that $\varphi(x)=\gamma(x)=1$ if $p=1$ or if $p=q$.

\begin{rem}
Formula \eqref{symform} should be compared with the property of locality in
vertex operator algebras \cite{B,FLM} (see also \cite{FF:laws}). Recall that
two power series $A(z)$ and $B(w)$ are called local if $A(z) B(w)$
converges for $|z|\gg|w|$, $B(w) A(z)$ converges for $|w|\gg|z|$, and their
analytic continuations satisfy: $A(z) B(w) = B(w) A(z)$.\qed
\end{rem}

The function $\varphi(x)$ is an elliptic function, i.e.
\begin{equation}    \label{fe1}
\varphi(xp^2) = \varphi(x),
\end{equation}
and it satisfies the functional equation
\begin{equation}    \label{fe2}
\varphi(x) \varphi(xp) = 1.
\end{equation}

The equations \eqref{fe1}, \eqref{fe2} provide a new understanding for
the algebra $\W_{q,p}(\sw_2)$. Let us explain that.

According to \eqref{conti1} the series
$\La(z)$ satisfies the relations:
\begin{equation}    \label{larel}
f_{1,2}\left(\frac{w}{z}\right) \La(z) \La(w) =
f_{1,2}\left(\frac{z}{w}\right) \La(w) \La(z).
\end{equation}

There is a difference between relations \eqref{equality} and
\eqref{larel}. The first is a relation on analytic continuations of the
compositions of two operators, while the second is a relation on formal
power series. A relation of the second type implies a relation of the first
type -- it can be obtained by multiplying it by a suitable meromorphic
function. But different relations of the second type may give rise to the
same relation of the first type as we will see below.

\begin{rem} Similar phenomenon occurs in vertex operator algebras. Consider
for example the Heisenberg algebra with generators $\beta_n,\gamma_n,
n\in\Z$, and relations $$[\beta_n,\gamma_m] = \kappa \delta_{n,-m},$$ where
$\kappa \in \C$. These relations imply the following formal power series
relations:
$$\beta(z) \gamma(w) - \gamma(w) \beta(z) = \kappa \delta\left( \frac{w}{z}
\right).$$ But we can also write
$$(z-w) \beta(z) \gamma(w) = (z-w) \gamma(w) \beta(z)$$ regardless of the
value of $\kappa$.\qed
\end{rem}

In order to get a relation of the second type from the relation
\eqref{equality} of the first type, we have to ``factorize'' the function
$\varphi(x)$, i.e. to represent it as $\varphi(x) = g(x^{-1})/g(x)$, where
$g(x)$ is a formal Taylor power series in $x$. Then we obtain a relation of
the second type $$g \left( \frac{w}{z} \right) \La(z) \La(w) = g \left(
\frac{z}{w} \right) \La(w) \La(z)$$ as in \eqref{larel}.

Such a factorization of $\varphi(x)$ is not unique in general. In our case
we can choose $g(x) = f_{1,2}(x)$, but we can also choose $g(x) =
f_{1,2}(xp^{-1})^{-1}$ or $g(x) = f_{1,2}(xp)^{-1}$, by virtue of the
equations \eqref{fe1} and \eqref{fe2}.

Let $\La_1(z)$ and $\La_2(z)$ satisfy the relations
\begin{align}    \label{fine1}
f_{1,2}\left(\frac{w}{z}\right) \La_i(z) \La_i(w) &=
f_{1,2}\left(\frac{z}{w}\right) \La_i(w) \La_i(z),\\    \label{fine2}
f_{1,2}\left(\frac{w}{zp}\right)^{-1} \La_1(z) \La_2(w) &=
f_{1,2}\left(\frac{zp}{w}\right)^{-1} \La_2(w) \La_1(z),\\    \label{fine3}
f_{1,2}\left(\frac{wp}{z}\right)^{-1} \La_2(z) \La_1(w) &=
f_{1,2}\left(\frac{z}{wp}\right)^{-1} \La_1(w) \La_2(z).
\end{align}

Then
\begin{equation}    \label{coarse}
\La_i(z) \La_j(w) = \varphi\left( \frac{w}{z} \right) \La_j(w)
\La_i(z)
\end{equation}
for all $i,j, \in \{ 1,2 \}$, i.e. elliptic relations \eqref{coarse} do not
depend on $i$ and $j$ while the formal power series relations
\eqref{fine1}--\eqref{fine3} do.

In fact, the relations \eqref{coarse} allow for even more sophisticated
formal power series relations. Let us set $T_1(z) = \La_1(z) +
\La_2(z)$. Then according to \eqref{coarse}, $T_1(z)$ satisfies the same
elliptic relations as $\La_i(z)$:
$$T_1(z) T_1(w) = \varphi\left( \frac{w}{z} \right) T_1(w) T_1(z).$$ But
when we write the relations between $T_1(z)$ and $T_1(w)$ as formal power
series, we obtain extra terms with $\delta$--functions:
\begin{multline}    \label{central}
f_{1,2} \left(\frac{w}{z}\right) T_1(z) T_1(w) - f_{1,2}
\left(\frac{z}{w}\right) T_1(w) T_1(z) \\ = \frac{(1-q)(1-p/q)}{1-p} \left(
\delta \left( \dfrac{w}{zp} \right) :\La_1(z) \La_2(zp): - \delta \left(
\dfrac{wp}{z} \right) :\La_1(w) \La_2(wp): \right).
\end{multline}

The reason is that the functions appearing in factorization of $\varphi(x)$
differ by certain rational multiples:
$$f_{1,2}(x) = f_{1,2}(xp)^{-1} \frac{(1-xq)(1-xp/q)}{(1-x)(1-xp)},$$
$$f_{1,2}(x) = f_{1,2}(xp^{-1})^{-1}
\frac{(1-xq^{-1})(1-xq/p)}{(1-x)(1-xp^{-1})}.$$ These multiples give rise
to the $\delta$--function terms in the right hand side of formula
\eqref{central}.

But we see from formula \eqref{fe2} that $:\La_1(z)\La_2(zp):$ is a central
element in the algebra generated by $\La_1(z)$ and $\La_2(z)$. Hence we may
set it equal to any number. If we set it equal $1$, then we obtain
relations \eqref{wqpsl2}, but we can put an arbitrary overall factor in the
right hand side of formula \eqref{wqpsl2}. In particular, if this factor is
0, we obtain the original defining relations \eqref{larel} of the
Heisenberg algebra. Hence the algebra $\W_{q,p}(\sw_2)$ is a central
extension of the Heisenberg algebra $\Hh_{q,p}(\sw_2)$.

The function given by \eqref{g} is the simplest solution of equations
\eqref{fe1} and \eqref{fe2}. It is interesting whether other solutions of
these equations give rise to quadratic algebras similar to
$\W_{q,p}(\sw_2)$.

\subsection{The case of $\sw_N$.}
{}From formulas \eqref{opel1}--\eqref{opel3} we find in the same way as in
the case of $\sw_2$:
\begin{equation}    \label{gN}
\La_i(z) \La_j(w) = \varphi_N\left(\frac{w}{z}\right) \La_j(w) \La_i(z),
\end{equation}
for all $i,j$, where $$\varphi_N(x) = \frac{\tN(x) \tN(xp) \tN(xp^{N-1}q)
\tN(xp^Nq^{-1})}{\tN(xq) \tN(xpq^{-1}) \tN(xp^{N-1}) \tN(xp^N)} =
\frac{\tN(xp) \tN(xq^{-1}) \tN(xp^{-1}q)}{\tN(xp^{-1}) \tN(xq)
\tN(xpq^{-1})}.$$ It is the last equality that ensures that the elliptic
relations between $\La_i(z)$ and $\La_j(w)$ are the same for $i=j$ and
$i\neq j$, although the ``factorized'' relations between them are
different, see formulas \eqref{opel1}--\eqref{opel3}.

Relations \eqref{gN} can be rewritten in a more symmetric form:
\begin{equation}    \label{symgN}
\gamma_N\left(\frac{w}{z}\right) \La_i(z) \La_j(w) =
\gamma_N\left(\frac{z}{w}\right) \La_j(w) \La_i(z),
\end{equation}
where $$\gamma_N(x) = \frac{\tN(x) \tN(xq)}{\tN(xp) \tN(xp^{-1}q)}.$$

Since relations \eqref{gN} do not depend on $i$ and $j$, we obtain similar
elliptic relations between for $T_i(z), i=1,\ldots,N-1$: $$T_i(z) T_j(w) =
\prod_{k=1}^i \prod_{l=1}^j \varphi_N\left( \frac{w}{z} p^{l-k} \right)
T_j(w) T_i(z).$$ But the formal power series relations between $T_i(z)$ and
$T_j(w)$ involve extra terms with $\delta$--functions as in the case of
$\sw_2$, see formula \eqref{m1}.

The function $\varphi_N(x)$ is elliptic, i.e. $\varphi_N(xp^N) =
\varphi_N(x)$, and it satisfies the functional equation
\begin{equation}    \label{feN}
\varphi_N(x) \varphi_N(xp) \ldots \varphi_N(xp^{N-1}) = 1.
\end{equation}
Equation \eqref{feN} implies that $:\La_1(z) \La_2(zp) \ldots
\La_N(zp^{N-1}):$ is a central element of $\Hh_{q,p}(\sw_N)$. Formula
\eqref{linrelp} shows that we have set it equal to $1$, but we could set it
equal to any number. If we do not set it equal to a number, we obtain
another algebra, which is natural to denote by $\W_{q,p}(\goth{gl}_N)$.

We see that in a certain sense the structure of the algebra
$\W_{q,p}(\sw_N)$ is ``encoded'' in equation \eqref{feN}, as in the case of
$\sw_2$.

\section{Elliptic relations for the screening currents.}
\label{elliptic}
\subsection{The screening currents of $\W_{q,p}(\sw_N)$.}
Let us we assume that $p$ and $q$ are generic and $|q|<1$. Then we have the
following relations for the screening currents when $|z|\gg|w|$:
\begin{align*}
S^+_i(z) S^+_i(w) & = z^{2\beta}
\frac{\left(\dfrac{w}{z}| 1,p;q\right)_\infty}{\left(\dfrac{w}{z}|
q,p^{-1}q;q\right)_\infty} :S^+_i(z) S^+_i(w):, \\
S^+_i(z) S^+_j(w) & = z^{-\beta}
\frac{\left(\dfrac{w}{z}| p^{-1/2}q;q\right)_\infty}{\left(\dfrac{w}{z}|
p^{1/2};q\right)_\infty} :S^+_i(z) S^+_j(w):, \quad \quad A_{ij}=-1,\\
S^+_i(z) S^+_j(w) & = :S^+_i(z) S^+_j(w):, \quad \quad A_{ij}=0.
\end{align*}

\begin{rem} If we set $p=q^{1-\beta}$, then in the limit $q \arr 1$ these
relations give us the well-known relations on the ordinary screening
currents:
$$S^+_i(z) S^+_j(w) = (z-w)^{A_{ij}\beta} :S^+_i(z) S^+_j(w):$$ when
$|z|\gg|w|$.\qed
\end{rem}

{}From the formulas above we obtain, in the analytic continuation sense,
\begin{align*}
S^+_i(z) S^+_i(w) & = - \left( \frac{w}{z} \right)^{1-2\beta}
\frac{\theta_q\left(\dfrac{w}{z}
p\right)}{\theta_q\left(\dfrac{z}{w}p\right)} S^+_i(w) S^+_i(z),\\ S^+_i(z)
S^+_j(w) & = \left( \frac{w}{z} \right)^{-2+\beta}
\frac{\theta_q\left(\dfrac{w}{z}
p^{-1/2}\right)}{\theta_q\left(\dfrac{z}{w}p^{-1/2}\right)} S^+_j(w)
S^+_i(z), \quad \quad A_{ij}=-1,\\
S^+_i(z) S^+_j(w) & = S^+_j(w) S^+_i(z), \quad \quad A_{ij}=0.
\end{align*}

We can rewrite these formulas as follows:
\begin{equation}    \label{final}
S^+_i(z) S^+_j(w) = (-1)^{A_{ij}-1} \left( \frac{w}{z}
\right)^{A_{ij}-A_{ij}\beta-1} \frac{\theta_q\left(\dfrac{w}{z}
p^{A_{ij}/2}\right)}{\theta_q\left(\dfrac{z}{w}p^{A_{ij}/2}\right)}
S^+_j(w) S^+_i(z).
\end{equation}

The function
\begin{equation}    \label{psi}
\psi_{ij}(x) = (-1)^{A_{ij}-1} x^{A_{ij}-A_{ij}\beta-1}
\frac{\theta_q\left( x p^{A_{ij}/2}\right)}{\theta_q\left(x^{-1}
p^{A_{ij}/2}\right)}
\end{equation}
appearing in the right hand side of formula \eqref{final} can be written as
\begin{equation}    \label{psidec}
\psi_{ij}(x) = \frac{\phi_{ij}(x)}{\phi_{ij}(x^{-1})},
\end{equation}
where
\begin{equation}    \label{phi}
\phi_{ij}(x) = x^{-\beta A_{ij}/2}
\frac{\theta_q(xp^{A_{ij}/2})}{\theta_q(xq^{A_{ij}/2})}.
\end{equation}

The functions $\psi_{ij}(x)$ and $\phi_{ij}(x)$ are quasi-periodic:
\begin{align*}
\psi_{ij}(xq) = \psi_{ij}(x), \quad \quad & \psi_{ij}(xe^{2\pi i}) =
e^{-2\pi i A_{ij}\beta} \psi_{ij}(x), \\ \phi_{ij}(xq) =
\phi_{ij}(x), \quad \quad & \phi_{ij}(xe^{2\pi i}) = e^{-\pi i A_{ij}
\beta} \phi_{ij}(x).
\end{align*}

Note also that $\psi_{ij}(x)=\phi_{ij}(x)=1$, if $p=q$.

To obtain relations on the screening currents $S^-_i(z), i=1,\ldots,N-1$,
let us assume that $p$ and $q$ are generic and $|p/q|<1$. Then we obtain in
the same way as above:
\begin{equation}    \label{final-}
S^-_i(z) S^-_j(w) = (-1)^{A_{ij}-1} \left( \frac{w}{z}
\right)^{A_{ij}-A_{ij}/\beta-1} \frac{\theta_{p/q}\left(\dfrac{w}{z}
p^{A_{ij}/2}\right)}{\theta_{p/q}\left(\dfrac{z}{w}p^{A_{ij}/2}\right)}
S^-_j(w) S^-_i(z).
\end{equation}

We also have:
\begin{align*}
S^+_i(z) S^-_i(w) & = \frac{1}{(z-wq)(z-wp^{-1}q)} :S^+_i(z) S^-_i(w):,\\
S^+_i(z) S^-_j(w) & = (z-wp^{-1/2}q) :S^+_i(z) S^-_j(w):, \quad \quad
A_{ij} = -1,\\ S^+_i(z) S^+_j(w) & = :S^+_i(z) S^-_j(w):, \quad \quad
A_{ij}=0.
\end{align*}

\subsection{The screening currents of $U_q(\sun)$.}
The screening currents involved in the Wakimoto realization of
$U_q(\sun)$ \cite{AOS} also satisfy elliptic relations.

These screening currents $S_i(z), i=1,\ldots,N-1$, are given in \cite{AOS}
by the formula
$$S_i(z) = :\exp \left( \sum_{n\neq 0} \frac{a^i_n}{[\nu n]_q} q^{-\nu
|n|/2} z^{-n} - \frac{1}{\nu} \hat{q}^i - \frac{1}{\nu} \hat{p}^i \log z
\right): \widetilde{S}_i(z).$$ The following relations hold \cite{AOS}:
$$[a^i_n,a^j_m] = \frac{1}{n} [\nu n]_q [A_{ij}]_q \delta_{n,-m}, \quad
\quad [\hat{p}^i,\hat{q}^j] = \nu A_{ij},$$ $$\widetilde{S}_i(z)
\widetilde{S}_j(w) = \frac{q^{-A_{ij}}z-w}{z-wq^{-A_{ij}}}
\widetilde{S}_j(w) \widetilde{S}_i(z),$$ in the sense of analytic
continuation, and $a^i_n$ commute with $\widetilde{S}_j(z)$. Here $\nu$ is
$k+g$ in the notation of \cite{AOS}.

Let us assume that $q$ is generic and $|q^{2\nu}|<1$. Then the relations
above give us the following relations on $S_i(z)$'s:
\begin{equation}    \label{wakskr}
S_i(z) S_j(w) = - \left( \frac{w}{z} \right)^{-A_{ij}/\nu-1}
\frac{\theta_t\left( \dfrac{w}{z} q^{-A_{ij}} \right)}{\theta_t\left(
\dfrac{z}{w} q^{-A_{ij}} \right)} S_j(w) S_i(z),
\end{equation}
where $t=q^{2\nu}$. These relations should be understood in the analytic
continuation sense.

The relations \eqref{final} can be viewed as elliptic deformations of the
quantum Serre relations of Drinfeld \cite{D}:
$$E_i(z) E_j(w) = \frac{zq^{A_{ij}}-w}{z-wq^{A_{ij}}} E_j(w) E_i(z),$$
which can be easily obtained from \eqref{wakskr} in the limit $t \arr 0$
with fixed $q$.

The function $$\Psi_{ij}(x) = - x^{-A_{ij}/\nu-1} \frac{\theta_t\left( x
q^{-A_{ij}} \right)}{\theta_t\left( x^{-1} q^{-A_{ij}} \right)}$$ in the
right hand side of formula \eqref{wakskr} can be rewritten as
$$\Psi_{ij}(x) = \frac{\Phi_{ij}(x)}{\Phi_{ij}(x^{-1})},$$ where
$$\Phi_{ij}(x) = x^{-A_{ij}/2\nu} \frac{\theta_t(x
q^{-A_{ij}})}{\theta_t(x)}.$$ The functions $\Psi_{ij}(x)$ and
$\Phi_{ij}(x)$ are quasi-periodic, and $\Psi_{ij}(x) = \Phi_{ij}(x) = 1$ in
the limit when $q \arr 1$ and $t$ is fixed.

\subsection{General case.}
Let $\g$ be a simply-laced simple Lie algebra. Let $\Hh'_{q,p}(\g)$ be the
Heisenberg algebra with generators $a_i[n], i=1,\ldots,l=\on{rank} \g;
n\in\Z$, and relations \eqref{ai}, where $(A_{ij})$ is the Cartan matrix of
$\g$. We define the Fock representations $\pi_\mu$ and the completion
$\Hh_{q,p}(\g)$ of $\Hh'_{q,p}(\g)$ in the same way as in
\secref{heisenberg}.

We define the screening currents $S^\pm_i(z), i=1,\ldots,l$, by formulas
\eqref{s+}--\eqref{sm-}. We then define the algebra $\W_{q,p}(\g)$ as the
subalgebra of $\Hh_{q,p}(\g)$ of elements which commute with the screening
currents $S^+_i(w)$ up to a total ${\cal D}_q$--difference. It follows from
the definition that $\W_{q,p}(\g)$ commutes with $S^-_i(w)$ up to a total
${\cal D}_{p/q}$--difference.

The relations between $S^+_i(z), i=1,\ldots,l$, in the analytic
continuation sense, are given by formula \eqref{final}. The relations
between $S^-_i(z), i=1,\ldots,l$, in the analytic continuation sense, are
given by formula \eqref{final-}. Note that formulas
\eqref{final}--\eqref{final-} make sense for an arbitrary integral
symmetric matrix $(A_{ij})$.

If we set $p=q^{1-\beta}$, then in the limit $q \arr 1$ the algebra
$\W_{q,p}(\g)$ becomes the ordinary $\W$--algebra, $\W_{\sqrt{\beta}}(\g)$
in the notation of \cite{FF:laws}. On the other hand, we expect that in
the limit $\beta \arr 0$, i.e. $p \arr q$, the algebra $\W_{q,p}(\g)$
becomes isomorphic to the Poisson algebra $\W_q(\g)$ considered in
\cite{FR}.

\vspace{10mm}
\noindent{\bf Acknowledgements.} We would like to thank T.~Miwa for the
invitation to visit R.I.M.S., and H.~Awata and A.~Odesskii for useful
discussions.

\end{document}